\documentclass[12pt,oneside]{amsart}
\usepackage{amssymb,verbatim,enumerate,amsmath,ifthen,cite}
\usepackage[mathscr]{eucal}
\usepackage[utf8]{inputenc}
\usepackage[T1]{fontenc}
\usepackage[marginparwidth=2.4cm]{geometry}
\newtheorem{theorem}{Theorem}[section]
\newtheorem*{theorem*}{Theorem}
\def\Thm#1#2{\ifthenelse{\equal{#1}{*}}{\begin{theorem*}#2\end{theorem*}}
             {\begin{theorem}\label{T#1}#2\end{theorem}}}
\newtheorem{Atheorem}{Theorem}

\def\thm#1{Theorem~\ref{T#1}}

\newtheorem{proposition}[theorem]{Proposition}
\newtheorem*{proposition*}{Proposition}
\def\Prp#1#2{\ifthenelse{\equal{#1}{*}}{\begin{proposition*}#2\end{proposition*}}
             {\begin{proposition}\label{P#1}#2\end{proposition}}}

\newtheorem{corollary}[theorem]{Corollary}
\newtheorem*{corollary*}{Corollary}
\def\Cor#1#2{\ifthenelse{\equal{#1}{*}}{\begin{corollary*}#2\end{corollary*}}
             {\begin{corollary}\label{C#1}#2\end{corollary}}}

\newtheorem{lemma}[theorem]{Lemma}
\newtheorem*{lemma*}{Lemma}
\def\Lem#1#2{\ifthenelse{\equal{#1}{*}}{\begin{lemma*}#2\end{lemma*}}
             {\begin{lemma}\label{L#1}#2\end{lemma}}}
\def\lem#1{Lemma~\ref{L#1}}
\newtheorem{Alemma}{Lemma}

\theoremstyle{definition}
\newtheorem{remark}[theorem]{Remark}
\newtheorem*{remark*}{Remark}
\def\Rem#1#2{\ifthenelse{\equal{#1}{*}}{\begin{remark*}\rm #2\end{remark*}}
             {\begin{remark}\label{R#1}\rm #2\end{remark}}}

\newtheorem{example}[theorem]{Example}
\newtheorem*{example*}{Example}
\def\Exa#1#2{\ifthenelse{\equal{#1}{*}}{\begin{example*}\rm #2\end{example*}}
             {\begin{example}\label{Ex#1}\rm #2\end{example}}}

\numberwithin{equation}{section}
\def\eq#1{{\rm(\ref{#1})}}
\def\Eq#1#2{\ifthenelse{\equal{#1}{*}}
  {\begin{equation*}\begin{aligned}[]#2\end{aligned}\end{equation*}}
  {\begin{equation}\begin{aligned}[]\label{#1}#2\end{aligned}\end{equation}}}

\newcommand{\operator}[1]{\mathop{\vphantom{\sum}\mathchoice
{\vcenter{\hbox{\LARGE $#1$}}}
{\vcenter{\hbox{\large $#1$}}}{#1}{#1}}\displaylimits}

\def\Mst_#1^#2{\operator{\mathscr{M}_{\mbox{\scriptsize$\#$}}\!\!}_{#1}^{#2}\,\,}

\newcommand\R{\mathbb{R}}

\newcommand\N{\mathbb{N}}

\DeclareMathOperator{\conv}{conv}

\def\comment#1{}

\newcounter{LT}

\title[Hölder- and Minkowski-type inequalities]{Hölder- and Minkowski-type inequalities \\ for generalized quasi-arithmetic means}

\author{Zsolt P\'ales}
\address{Institute of Mathematics, University of Debrecen, 
Pf.\ 400, 4002 Debrecen, Hungary}
\email{pales@science.unideb.hu}

\author{Pawe{\l} Pasteczka}
\address{Institute of Mathematics, University of Rzesz\'ow, Pigonia 1, 35-310 Rzesz\'ow, Poland}
\email{ppasteczka@ur.edu.pl}



\begin{document}

\begin{abstract}
The purpose of this paper is to establish several necessary and sufficient conditions to ensure the validity of a general functional inequality in terms of generalized quasi-arithmetic means. In particular cases, we consider  
Hölder-, Minkowski-, and Jensen-type inequalities. Generalized quasi-arithmetic means are defined by taking strictly monotone generating functions instead of strictly monotone and continuous ones. 
\end{abstract}
\maketitle

\section{Introduction}

This paper aims to study the functional inequality of the form 
\Eq{GFI}{  M_0\big(\varphi(x_{1,1},\dots,x_{1,k}),\dots,\varphi(x_{n,1},\dots,x_{n,k})\big)
  \leq \Phi\Big(M_1(x_{1,1},\dots,x_{n,1}),\dots,M_k(x_{1,k},\dots,x_{n,k})\Big),\\
  \text{ for all }n\in\N\text{ and for all }
  (x_{1,1},\dots,x_{1,k}),\dots,(x_{n,1},\dots,x_{n,k})\in I_1\times\dots\times I_k,
}
where $n,k\in \N$, $M_j\colon I_j^n \to I_j$ are $n$-variable means on $I_j$ ($j \in \{0,\dots,k\}$), $\varphi\colon I_1\times\dots\times I_k \to I_0$, $\Phi \colon I_1\times\dots\times I_k \to \R$ ($I_0,\dots,I_k$ are subintervals of $\R$) in the case when all $M_i$-s are generalized quasi-arithmetic means. 
More precisely, we characterize, for a given pair $(\varphi,\Phi)$, all sequences of generalized quasi-arithmetic means $(M_0,\dots,M_k)$ satisfying this inequality.

Furthermore, we will study its weighted counterpart, that is we involve the space of weights $\Lambda_n:=[0,\infty)^n \setminus\{(0,\dots,0)\}$ in the inequality \eq{GFI} and we obtain
\Eq{GwFI}{
  \widetilde{M}_0&\big(\varphi(x_{1,1},\dots,x_{1,k}),\dots,\varphi(x_{n,1},\dots,x_{n,k});\lambda_1,\dots,\lambda_n\big)\\
  &\hspace{5mm}\leq \Phi\Big(\widetilde{M}_1(x_{1,1},\dots,x_{n,1};\lambda_1,\dots,\lambda_n),\,\dots,\,\widetilde{M}_k(x_{1,k},\dots,x_{n,k};\lambda_1,\dots,\lambda_n)\Big),\\
  &\hspace{10mm}\text{ for all }n\in\N,\ \lambda \in W_n, \text{ and } (x_{1,1},\dots,x_{1,k}),\dots,(x_{n,1},\dots,x_{n,k})\in I_1\times\dots\times I_k,
}
where all remaining parameters are the same as above except $\widetilde{M}_j \colon I_j^n \times \Lambda_n \to I_j$ ($j\in\{0,\dots,k\}$) which, in this case, are $n$-variable weighted means on suitable intervals.

Inequality \eq{GwFI} seems to be more general than \eq{GFI} as it contains \eq{GFI} as a special case (when $\lambda\equiv 1)$. Nevertheless, these two approaches are equivalent in view of our main result. These inequalities generalize several known problems. For example, if $\varphi=\Phi$ is the arithmetic mean, then we get the convexity-type inequalities; if $\varphi=\Phi$ is the sum, then our problem reduces to Minkowski-type inequalities, and so on. It will turn out that $\varphi=\Phi$ is generally the important subclass of these inequalities providing meaningful simplifications.

Let us also emphasize that converse inequalities can be obtained using the standard reflection-type techniques. Roughly speaking, we consider reflected means $M_j^*$ defined on $-I_j$ as $M_j^*(x):=-M_j(-x)$ instead of the original ones and apply natural substitutions; see for example \cite[sect.~2.2]{PalPas21a}. We skip the details for the sake of compactness.

Finally, we refer to the papers \cite{GruPal24,Bec70,CziPal00,CziPal03,LosPal97a,LosPal98,Pal82b,Pal83c,Pal86a,Pal88a,Pal99a}, where the reader can find results related to the standard quasi-arithmetic setting.

\section{The concept of generalized inverse and generalized quasi-arithmetic means}
The concept of generalized quasi-arithmetic means was mentioned for the first time in \cite[Example 4.5(3)]{PalPas19a} as an important subclass of semideviation means (a large class of means defined for the first time in \cite{Pal89b}). Furthermore, these means are the particular cases of the generalized Bajraktarević means introduced in \cite{GruPal20}. The first comprehensive study of this family was presented in the recent paper \cite{PalPas2412}.


Here we recall some relevant concepts and results from the past. In the sequel, the set of points in $I$, where $f$ is continuous, will be denoted by $C_f$.

\Lem{SMF}{{\rm(\!\!\cite[Lemma 1]{GruPal20})} Let $f:I\to\R$ be a strictly increasing function. Then there exists a uniquely determined increasing function $f^{(-1)}:\conv(f(I))\to I$ which is the left inverse of $f$, i.e., for all $x\in I$,
\Eq{SMF}{
   (f^{(-1)}\circ f)(x)=x.
}
Furthermore, 
$f^{(-1)}$ is continuous and $(f\circ f^{(-1)})(y)=y$ for all $y\in f(I)$.
}

This formulation of the latter statement follows the one presented in  \cite{PalPas2412}. The function $f^{(-1)}$ described in the above lemma is called the \emph{generalized inverse of the strictly increasing function $f:I\to\R$}. It is clear that the restriction of $f^{(-1)}$ to $f(I)$ equals the inverse of $f$ in the standard sense. Therefore, $f^{(-1)}$ is the continuous and increasing extension of the inverse of $f$ to the smallest interval containing the range of $f$. It turns out that the generalized inverse follows the composition principle known for the common inverse. This property, however, requires a simple proof.
\Lem{GI}{
Let $I,J$ be nonempty open subintervals of $\R$, and let $f\colon I \to J$, $g \colon J \to \R$ be strictly increasing functions. Then $(f\circ g)^{(-1)}=g^{(-1)}\circ f^{(-1)}$.
}
\begin{proof}
Note that $f\circ g$ is strictly increasing.  Take any $y_0 \in \conv(f \circ g(I))$, and set $x_0:=(f\circ g)^{-1}(y_0)$. Then for all $x,x' \in I$ with $x<x_0<x'$ we have
$f\circ g(x) < y_0 < f\circ g(x')$.

Applying, side-by-side, first $f^{(-1)}$,  then $g^{(-1)}$, and using that these functions are increasing, we get $x \le g^{(-1)} \circ f^{(-1)} (y_0) \le x'$. Finally, upon taking the limit $x \uparrow x_0$ and $x' \downarrow x_0$ we get 
$g^{(-1)} \circ f^{(-1)} (y_0)=x_0=(f\circ g)^{-1}(y_0)$,
which completes the proof.
\end{proof}

Let us now recall the easy but useful lemma.

\Lem{ineq}{{\rm(\!\!\cite[Lemma 2.2]{PalPas2412})} Let $f:I\to\R$ be a strictly increasing function and let $x\in I$ and $u\in\conv(f(I))$. Then we have the following equivalences:
\begin{enumerate}[(i)]
 \item $f^{(-1)}(u)<x$ holds if and only if $u<f_-(x)$.
 \item $f^{(-1)}(u)\leq x$ holds if and only if $u\leq f_+(x)$.
 \item $x<f^{(-1)}(u)$ holds if and only if $f_+(x)<u$.
 \item $x\leq f^{(-1)}(u)$ holds if and only if $f_-(x)\leq u$.
\end{enumerate}
}

Given a strictly increasing function $f:I\to\R$ and $n\in\N$, the \emph{$n$-variable generalized quasi-arithmetic mean} $A^{[n]}_{f}:I^n\to I$ and the \emph{$n$-variable weighted generalized quasi-arithmetic mean} $\widetilde{A}^{[n]}_{f}:I^n\times\Lambda_n\to I$ are defined by the following formulas:
\Eq{QAM}{
	A^{[n]}_{f}(x):=f^{(-1)}\bigg(\frac{f(x_1)+\dots+f(x_n)}{n}\bigg) \qquad (x\in I^n),
}
\Eq{WQAM}{
	\widetilde{A}^{[n]}_{f}(x;\lambda):=f^{(-1)}\bigg(\frac{\lambda_1f(x_1)+\dots+\lambda_nf(x_n)}{\lambda_1+\dots+\lambda_n}\bigg) \qquad (x\in I^n,\,\lambda\in\Lambda_n).
}
Clearly, for any $x\in I^n$, we have that $A^{[n]}_{f}(x)=\widetilde{A}^{[n]}_{f}(x;\lambda)$ if $\lambda_1=\dots=\lambda_n$.

\Lem{SC}{{\rm(\!\!\cite[Proposition 3.3]{PalPas2412})}
Let $n\in\N\setminus\{1\}$, let $f,g:I\to\R$ be strictly increasing functions, and let $t\in I$. Assume that $A_f^{[n]}(x)\leq A_g^{[n]}(x)$ holds for all $x\in I^n$. Then, we have the following assertions:
\begin{enumerate}[(i)]
 \item If $g$ is lower semicontinuous at $t$, then so is $f$.
 \item If $f$ is upper semicontinuous at $t$, then so is $g$.
 \item The function $f$ is continuous at $t$ if and only if $g$ is continuous at $t$.
\end{enumerate}}

\section{Main results}

Our main result reads as follows.

\Thm{Main}{Let $k\in\N$, let $I_0,I_1,\dots,I_k$ be nonempty open intervals and denote $I:=I_1\times\dots\times I_k$. For all $j\in\{0,1,\dots,k\}$, let $f_j:I_j\to\R$ be strictly increasing, let $\varphi:I\to I_0$ be arbitrary and let $\Phi:I\to\R$ be separately increasing and upper semicontinuous. Assume that the set $\Gamma:=\Phi^{-1}(C_{f_0})=\{t\in I\colon \Phi(t)\in C_{f_0}\}$ is a dense subset of $I$. Then the following assertions are equivalent.
\begin{enumerate}[(a)]
\item For all $n\in\N$ and for all $(x_{1,1},\dots,x_{1,k}),\dots,(x_{n,1},\dots,x_{n,k})\in I$, 
\Eq{*}{
  A^{[n]}_{f_0}\big(\varphi(x_{1,1},\dots,x_{1,k}),\dots,\varphi(x_{n,1},\dots,x_{n,k})\big)
  \leq \Phi\Big(A^{[n]}_{f_1}(x_{1,1},\dots,x_{n,1}),\dots,A^{[n]}_{f_k}(x_{1,k},\dots,x_{n,k})\Big).
}
\item For all $n\in\N$, for all $(x_{1,1},\dots,x_{1,k}),\dots,(x_{n,1},\dots,x_{n,k})\in I$, and for all $(\lambda_1,\dots,\lambda_n)\in\Lambda_n$, 
\Eq{*}{
  \widetilde{A}^{[n]}_{f_0}\big(\varphi(x_{1,1},\dots,x_{1,k})&,\dots,\varphi(x_{n,1},\dots,x_{n,k});\lambda_1,\dots,\lambda_n\big)\\
  &\leq \Phi\Big(\widetilde{A}^{[n]}_{f_1}(x_{1,1},\dots,x_{n,1};\lambda_1,\dots,\lambda_n),\,\dots,\,\widetilde{A}^{[n]}_{f_k}(x_{1,k},\dots,x_{n,k};\lambda_1,\dots,\lambda_n)\Big).
}
\item For all $(x_{1,1},\dots,x_{1,k}),\dots,(x_{k+1,1},\dots,x_{k+1,k})\in I$ and for all $(\lambda_1,\dots,\lambda_{k+1})\in\Lambda_{k+1}$, 
\Eq{*}{
  \widetilde{A}^{[k+1]}_{f_0}&\big(\varphi(x_{1,1},\dots,x_{1,k}),\dots,\varphi(x_{k+1,1},\dots,x_{k+1,k});\lambda_1,\dots,\lambda_{k+1}\big)\\
  &\leq \Phi\Big(\widetilde{A}^{[k+1]}_{f_1}(x_{1,1},\dots,x_{k+1,1};\lambda_1,\dots,\lambda_{k+1}),\,\dots,\,\widetilde{A}^{[k+1]}_{f_k}(x_{1,k},\dots,x_{k+1,k};\lambda_1,\dots,\lambda_{k+1})\Big).
}
\item For all $j\in\{1,\dots,k\}$, there exists a nonnegative function $a_j:\Gamma\to\R$ such that
for all $(x_1,\dots,x_k)\in I$ and $(t_1,\dots,t_k)\in \Gamma$,
\Eq{*}{
f_0\big(\varphi(x_1,\dots,x_k)\big)-f_0\big(\Phi(t_1,\dots,t_k)\big)\leq \sum_{j=1}^k a_j(t_1,\dots,t_k)\big(f_j(x_j)-f_j(t_j)\big).
}
\item There exists a separately increasing concave function $\Psi:\conv(f_1(I_1))\times\dots\times\conv(f_k(I_k))\to \R$ such that, for all $(x_1,\dots,x_k)\in I$,
\Eq{e1}{
  f_0\big(\varphi(x_1,\dots,x_k)\big)
  \leq\Psi\big(f_1(x_1),\dots,f_k(x_k)\big)
}
and, for all $(t_1,\dots,t_k)\in\Gamma$,
\Eq{e2}{
  f_0\big(\Phi(t_1,\dots,t_k)\big)
  \geq\Psi\big(f_1(t_1),\dots,f_k(t_k)\big).
}
\end{enumerate}
}

\begin{proof} The proof of implication $(a)\Rightarrow(b)$ is standard, it is based on the continuity of generalized quasi-arithmetic means with respect to weights and the upper semicontinuity of $\Phi$.

Applying assertion $(b)$ for $n=k+1$, we can see that $(b)$ implies $(c)$.

To verify the implication $(c)\Rightarrow(d)$, 
assume that the assertion $(c)$ is valid. To prove the existence of nonnegative functions described in assertion $(d)$ it is sufficient to show that, for all fixed $t:=(t_1,\dots,t_k)\in \Gamma$, there exists $(a_1,\dots,a_k)\in[0,\infty)^k$ such that
\Eq{*}{
f_0\big(\varphi(x)\big)-f_0\big(\Phi(t)\big)\leq \sum_{j=1}^k a_j\big(f_j(x_j)-f_j(t_j)\big)
}
holds for all $x=(x_1,\dots,x_k)\in I$. To see this, we prove that the intersection of the family of closed convex sets $\{H(x)\colon x\in I\}$ is not empty, where
\Eq{*}{
  H(x):=\bigg\{(a_1,\dots,a_k)\in[0,\infty)^k\colon f_0\big(\varphi(x)\big)-f_0\big(\Phi(t)\big)\leq \sum_{j=1}^k a_j\big(f_j(x_j)-f_j(t_j)\big)\bigg\},\qquad x\in I.
}
First we show that there exists a finite (in fact a $k$-element) subset of $I$ over which the intersection is bounded and hence it is compact.

Let $e_1,\dots,e_k$ denote the standard basis for $\R^k$. Choose $\varepsilon>0$ so that $t-\varepsilon e_j\in I$ be valid for all $j\in\{1,\dots,k\}$. Then
\Eq{*}{
  H(t-\varepsilon e_j)=\bigg\{(a_1,\dots,a_k)\in[0,\infty)^k\colon f_0\big(\varphi(t-\varepsilon e_j)\big)-f_0\big(\Phi(t)\big)\leq a_j\big(f_j(t_j-\varepsilon)-f_j(t_j)\big)\bigg\},
}
which, by the strict monotonicity of $f_j$, shows that, for $(a_1,\dots,a_k)\in H(t-\varepsilon e_j)$, we have
\Eq{*}{
   a_j\in\bigg[0,\frac{f_0\big(\varphi(t-\varepsilon e_j)\big)-f_0\big(\Phi(t)\big)}{f_j(t_j-\varepsilon)-f_j(t_j)}\bigg].
}
Therefore,
\Eq{*}{
  \bigcap_{j=1}^k H(t-\varepsilon e_j)
  \subseteq \prod_{j=1}^k \bigg[0,\frac{f_0\big(\varphi(t-\varepsilon e_j)\big)-f_0\big(\Phi(t)\big)}{f_j(t_j-\varepsilon)-f_j(t_j)}\bigg].
}
According to an infinite variant of Helly's Theorem, if the intersection $\bigcap_{x\in I}H(x)$ were empty, then would there exists a $(k+1)$-member subfamily whose intersection is also empty. That is, there exist
$(x_{1,1},\dots,x_{1,k})$, \dots, $(x_{k+1,1},\dots,x_{k+1,k})\in I$ such that the following system of linear inequalities
\Eq{*}{
f_0\big(\varphi(x_{i,1},\dots,x_{i,k})\big)-f_0\big(\Phi(t)\big)\leq \sum_{j=1}^k a_j\big(f_j(x_{i,j})-f_j(t_j)\big)
\qquad (i\in\{1,\dots,k+1\})
}
has no solution for $(a_1,\dots,a_k)\in[0,\infty)^k$. Using a standard version of the celebrated Farkas Lemma, it follows that there exists $\lambda=(\lambda_1,\dots,\lambda_{k+1})\in\Lambda_{k+1}$ such that
\begin{align}
  &\sum_{i=1}^{k+1} \lambda_i \big(f_j(x_{i,j})-f_j(t_j)\big)\leq 0 \qquad(j\in\{1,\dots,k\})
  \qquad\mbox{and}\label{I:j}\\
  &\sum_{i=1}^{k+1} \lambda_i\big(f_0\big(\varphi(x_{i,1},\dots,x_{i,k})\big)-f_0\big(\Phi(t)\big)\big)>0. \label{I:0}
\end{align}
For $j\in\{1,\dots,k\}$, the inequality \eqref{I:j} implies that
\Eq{*}{
  \frac{\sum_{i=1}^{k+1} \lambda_i f_j(x_{i,j})}{\sum_{i=1}^{k+1} \lambda_i}\leq f_j(t_j),
}
which then yields that
\Eq{*}{
  \widetilde{A}^{[k+1]}_{f_j}\big(x_{1,j},\dots,x_{k+1,j};\lambda_1,\dots,\lambda_{k+1}\big)=f_j^{(-1)}\bigg(\frac{\sum_{i=1}^{k+1} \lambda_i f_j(x_{i,j})}{\sum_{i=1}^{k+1} \lambda_i}\bigg)\leq f_j^{(-1)}(f_j(t_j))=t_j.
}
From the inequality \eqref{I:0} it follows that
\Eq{*}{
  \frac{\sum_{i=1}^{k+1} \lambda_if_0\big(\varphi(x_{i,1},\dots,x_{i,k})\big)}{\sum_{i=1}^{k+1} \lambda_i}> f_0\big(\Phi(t)\big).
}
Therefore, using that $\Phi(t)\in C_{f_0}$, we get
\Eq{*}{
  &\widetilde{A}^{[k+1]}_{f_0}\big(\varphi(x_{1,1},\dots,x_{1,k}),\dots ,\varphi(x_{k+1,1},\dots,x_{k+1,k});\lambda_1,\dots,\lambda_{k+1}\big)\\
  &=f_0^{(-1)}\bigg(\frac{\sum_{i=1}^{k+1} \lambda_if_0\big(\varphi(x_{i,1},\dots,x_{i,k})\big)}{\sum_{i=1}^{k+1} \lambda_i}\bigg)>\Phi(t)\\
  &\geq\Phi\Big(\widetilde{A}^{[k+1]}_{f_1}(x_{1,1},\dots,x_{k+1,1};\lambda_1,\dots,\lambda_{k+1}),\,\dots,\,\widetilde{A}^{[k+1]}_{f_k}(x_{1,k},\dots,x_{k+1,k};\lambda_1,\dots,\lambda_{k+1})\Big),
}
which contradicts assumption $(c)$. The contradiction so obtained shows that the intersection $\bigcap_{x\in I}H(x)$ is not empty.

In the next step, we show that assertion $(d)$ implies $(e)$.
Assume that $(d)$ holds and, for $(t_1,\dots,t_k)\in \Gamma$, define the affine function $\Psi_{t_1,\dots,t_k}:\R^k\to\R$ by
\Eq{*}{
  \Psi_{t_1,\dots,t_k}(u_1,\dots,u_k)
  :=f_0\big(\Phi(t_1,\dots,t_k)\big)+\sum_{j=1}^k a_j(t_1,\dots,t_k)\big(u_j-f_j(t_j)\big).
}
Due to the nonnegativity of $a_1(t_1,\dots,t_k),\dots,a_k(t_1,\dots,t_k)$, it follows that $\Psi_{t_1,\dots,t_k}$ is separately increasing.
According to condition $(d)$, 
for all $(x_1,\dots,x_k)\in I$, we also have that
\Eq{e1+}{
  f_0\big(\varphi(x_1,\dots,x_k)\big)
  \leq\Psi_{t_1,\dots,t_k}\big(f_1(x_1),\dots,f_k(x_k)\big)
}
and
\Eq{e2+}{
  f_0\big(\Phi(t_1,\dots,t_k)\big)
  =\Psi_{t_1,\dots,t_k}\big(f_1(t_1),\dots,f_k(t_k)\big).
}
Now we define $\Psi:\conv(f_1(I_1))\times\dots \times\conv(f_k(I_k))\to[-\infty,\infty]$ as follows
\Eq{*}{
  \Psi(u_1,\dots,u_k)
  :=\inf_{(t_1,\dots,t_k)\in \Gamma} \Psi_{t_1,\dots,t_k}(u_1,\dots,u_k).
}
Then $\Psi$ is concave and separately increasing because it is the pointwise infimum of affine (and hence also concave) separately increasing functions. 
For any fixed $(x_1,\dots,x_k)\in I$, taking the infimum of the right-hand side in inequality \eq{e1+} with respect to $(t_1,\dots,t_k)\in \Gamma$, we can see that \eq{e1} holds.
It is obvious that $\Psi(u)<\infty$ for all $u\in \conv(f_1(I_1))\times\dots\times\conv(f_k(I_k))$. By the inequality \eq{e1+}, we have that $\Psi(u)>-\infty$ for all $u\in f_1(I_1)\times\dots\times f_k(I_k)$. Using the concavity of $\Psi$, it follows that the inequality $\Psi(u)>-\infty$ is true for all $u$ belonging to the convex hull of the set $f_1(I_1)\times\dots\times f_k(I_k)$, which equals $\conv(f_1(I_1))\times\dots\times\conv(f_k(I_k))$. Hence, the function $\Psi$ takes only finite values over its domain. The equality \eq{e2+} and the inequality $\Psi_{t_1,\dots,t_k}\geq\Psi$ show that the inequality \eq{e2} is valid for all $(t_1,\dots,t_k)\in\Gamma$.

Finally, assume that $(e)$ holds. To show that $(a)$ is valid, let $n\in\N$, and let $(x_{i,1},\dots,x_{i,k})\in I$ for $i\in\{1,\dots,n\}$. Let $(t_1,\dots,t_k)\in\Gamma$ be arbitrary such that the inequality
\Eq{ejj}{
  A^{[n]}_{f_j}(x_{1,j},\dots,x_{n,j})<t_j
  \qquad(j\in\{1,\dots,k\})
}
be valid. This implies that
\Eq{ej}{
\frac{1}n\sum_{i=1}^n f_j(x_{i,j})< f_j(t_j)
\qquad(j\in\{1,\dots,k\}).
}
On the other hand, by inequality \eq{e1}, we have
\Eq{*}{
  f_0\big(\varphi(x_{i,1},\dots,x_{i,k})\big)
  \leq \Psi\big(f_1(x_{i,1}),\dots,f_k(x_{i,k})\big)
  \qquad i\in\{1,\dots,n\}.
}
Therefore, using the concavity first, then the separate increasingness of $\Psi$ together with the inequalities in \eq{ej}, finally applying the inequality \eq{e2}, we conclude that
\Eq{*}{
  \frac{1}{n}\sum_{i=1}^n
  f_0\big(\varphi(x_{i,1},\dots,x_{i,k})\big)
  &\leq \frac{1}{n}\sum_{i=1}^n
  \Psi\big(f_1(x_{i,1}),\dots,f_k(x_{i,k})\big)\\
  &\leq\Psi\bigg(\frac{1}{n}\sum_{i=1}^nf_1(x_{i,1}),\dots,\frac{1}{n}\sum_{i=1}^nf_k(x_{i,k})\bigg)\\
  &\leq\Psi\big(f_1(t_1),\dots,f_k(t_k)\big)
  \leq f_0\big(\Phi(t_1,\dots,t_k)\big).
}
Thus, we have proved that, for all $(t_1,\dots,t_k)\in\Gamma$ which satisfies the inequalities in \eq{ejj}, 
\Eq{*}{
  \frac{1}{n}\sum_{i=1}^n
  f_0\big(\varphi(x_{i,1},\dots,x_{i,k})\big)
  \leq f_0\big(\Phi(t_1,\dots,t_k)\big).
}
Hence, it follows that
\Eq{*}{
  A^{[n]}_{f_0}\big(\varphi(x_{1,1},\dots,x_{1,k}),\dots,\varphi(x_{n,1},\dots,x_{n,k})\big)
  =f_0^{(-1)}\bigg(\frac{1}{n}\sum_{i=1}^n
  f_0\big(\varphi(x_{i,1},\dots,x_{i,k})\big)\bigg)
  \leq \Phi(t_1,\dots,t_k).
}
Finally, we apply that the set $\Gamma$ is dense in $I$ and $\Phi$ is upper semicontinuous on this product set. Hence, we can take the upper limit as 
\Eq{*}{
\Gamma \ni (t_1,\dots,t_k)\to\Big(A^{[n]}_{f_1}(x_{1,1},\dots,x_{n,1}),\dots,A^{[n]}_{f_k}(x_{1,k},\dots,x_{n,k})\Big)
}
in the above inequality, we obtain that assertion $(a)$ is valid.
\end{proof}

\Rem{*}{If we assume that the function $\Phi:I\to\R$ is strictly increasing in one of its variables, then the set $\Gamma=\Phi^{-1}(C_{f_0})$ is dense in $I$. Indeed, assume that $\Phi$ is strictly increasing in its $k$th variable. Let $(t_1,\dots,t_{k-1})\in I_1\times\dots\times I_{k-1}$ be arbitrary. Then the map $t_k\mapsto \varphi(t_1,\dots,t_{k-1},t_k)$ is an injective mapping of $I_k$ into $I_0$. Therefore, the inverse image by this map of the set of discontinuity points of $f_0$ which is countable is also countable. Therefore, the inverse image by this map of the set of continuity points of $f_0$ is co-countable and hence it is dense in $I_k$. This proves that $\Gamma$ is dense in $I$.}

\Thm{cont}{
Let $k\in\N$, let $I_0,I_1,\dots,I_k$ be nonempty open intervals and denote $I:=I_1\times\dots\times I_k$. For all $j\in\{0,1,\dots,k\}$, let $f_j:I_j\to\R$ be strictly increasing and let $\varphi:I\to I_0$ be separately strictly increasing. Assume, for some $n\in\N\setminus\{1\}$, that for all $(x_{1,1},\dots,x_{1,k}),\dots,(x_{n,1},\dots,x_{n,k})\in I$, the inequality
\Eq{347}{
  A^{[n]}_{f_0}\big(\varphi(x_{1,1},\dots,x_{1,k}),\dots,\varphi(x_{n,1},\dots,x_{n,k})\big)
  \leq \varphi\Big(A^{[n]}_{f_1}(x_{1,1},\dots,x_{n,1}),\dots,A^{[n]}_{f_k}(x_{1,k},\dots,x_{n,k})\Big)
}
holds. Then, for all $j\in\{1,\dots,k\}$, the function $f_j$ is continuous on $I_j$ and $f_0$ is continuous on $\varphi(I)$. Furthermore, $\varphi$ is also continuous on $I$.
}

\begin{proof}
Fix $\ell \in \{1,\dots,k\}$ arbitrarily. We show that $f_{\ell} \colon I_{\ell} \to \R$ is continuous. For a given vector $z=(z_1,\dots,z_k)\in I$ let us define a function $\psi_{z,\ell} \colon I_{\ell} \to I_0$ by 
\Eq{*}{
\psi_{z,\ell}(x):=\varphi(z_1,\dots,z_{\ell-1},x,z_{\ell+1},\dots,z_k),
}
that is, $x$ appears at the $\ell$-th coordinate of the argument of $\varphi$. Since $\varphi$ is separately strictly increasing, the function $\psi_{z,\ell}$ is strictly increasing on $I_{\ell}$ (for all $z \in I)$. In particular, $\psi_{z,\ell}$ possesses an increasing left inverse $\psi_{z,\ell}^{(-1)}\colon \conv(\psi_{z,\ell}(I_\ell))\to I_\ell$. For an arbitrary $y \in I_{\ell}^n$ let us substitute 
\Eq{*}{
x_{i,j}=\begin{cases}
y_i &\quad \text{ for }(i,j)\in \{1,\dots,n\}\times \{\ell\}\\
z_j &\quad \text{ for }(i,j)\in \{1,\dots,n\} \times \big(\{1,\dots,k\}\setminus\{\ell\}\big)
\end{cases}
}
into \eq{347}. Then, for $j\in \{1,\dots,k\}\setminus\{\ell\}$, we have that $A^{[n]}_{f_j}(x_{1,j},\dots,x_{n,j})=A^{[n]}_{f_j}(z_j,\dots,z_j)=z_j$, hence \eq{347} simplifies to
\Eq{367}{
A_{f_0}^{[n]}(\psi_{z,\ell}(y_1),\dots,\psi_{z,\ell}(y_n))\le \psi_{z,\ell}(A_{f_{\ell}}^{[n]}(y_1,\dots,y_n)) \qquad (y_1,\dots,y_n)\in I_{\ell}^n.
}

Therefore, since $\psi_{z,\ell}^{(-1)}$ is an increasing function with $\psi_{z,\ell}^{(-1)}\circ \psi_{z,\ell}=\textrm{id}$, we get
\Eq{*}{
\psi_{z,\ell}^{(-1)}\circ A_{f_0}^{[n]}(\psi_{z,\ell}(y_1),\dots,\psi_{z,\ell}(y_n))\le A_{f_{\ell}}^{[n]}(y_1,\dots,y_n) \qquad (y_1,\dots,y_n)\in I_{\ell}^n.
}
In view of \lem{GI} we get know that  $\psi_{z,\ell}^{(-1)}\circ f_0^{(-1)}=(f_0\circ \psi_{z,\ell})^{(-1)}$, and therefore
\Eq{*}{
A_{f_0\circ \psi_{z,\ell}}^{[n]}(y_1,\dots,y_n)\le A_{f_{\ell}}^{[n]}(y_1,\dots,y_n) \qquad (y_1,\dots,y_n)\in I_{\ell}^n.
}
Therefore, in view of \lem{SC}, we get that $C_{f_0 \circ \psi_{z,\ell}}=C_{f_{\ell}}$ holds for all $z \in I$. 

Now assume, to the contrary, that $f_{\ell}$ is not continuous at a certain point $x_0 \in I_{\ell}$. Then $f_{0}$ is not continuous at all points belonging to the set
\Eq{*}{
S:=\{\psi_{z,\ell}(x_0)\colon z \in I\}=
\{\varphi(z_1,\dots,z_{\ell-1},x_0,z_{\ell+1},\dots,z_k) \colon z \in I\}.
}
Since $z\mapsto \psi_{z,\ell}(x_0)$ is separately strictly increasing on $I$ except with respect to the $\ell$th variable, we find that $S$ is an uncountable set, i.e. $f_0$ is not continuous at uncountably many points, which is impossible as $f_0$ was assumed to be strictly increasing. This contradiction shows that $f_{\ell}$ must be continuous on $I_{\ell}$ and hence $C_{f_{\ell}}=I_{\ell}$. Consequently, $C_{f_0 \circ \psi_{z,\ell}}=I_{\ell}$, which shows that $f_0 \circ \psi_{z,\ell}$ is continuous on $I_{\ell}$ for all $z\in I$. Having in mind that $f_0,\psi_{z,\ell}$ are strictly increasing on their domains, it follows that 
$\psi_{z,\ell}$ is continuous on $I_{\ell}$ and $f_0$ is continuous on the open interval $\psi_{z,\ell}(I_{\ell})$ for all $z\in I$. Therefore, $f_0$ is continuous over the union $$\bigcup_{z\in I}\psi_{z,\ell}(I_{\ell})=\varphi(I).$$

Finally, using that $\psi_{z,\ell}$ is continuous on $I_{\ell}$ for all $\ell\in\{1,\dots,k\}$ and $z\in I$, we show that $\varphi$ is continuous on $I$.
Let $(p_1,\dots,p_k)\in I$ be a point at which we want to verify the continuity of $\varphi$. Let $\varepsilon>0$ be arbitrary. Let $e_1,\dots,e_k$ be the standard basis in $\R^k$. We will recursively construct the positive numbers $\delta_1,\dots,\delta_k>0$ as follows: Let $\ell\in\{1,\dots,k\}$ and assume that $\delta_j$ has been constructed for $j\in\{1,\dots,\ell-1\}$. Then the map
\Eq{*}{
 I_\ell-p_\ell\ni x\mapsto \psi_{z,\ell}(p_\ell+x)=\varphi\bigg(p+\sum_{j=1}^{\ell-1}\delta_je_j+xe_\ell\bigg)
}
(where $z=p+\sum_{j=1}^{\ell-1}\delta_je_j$) is continuous and strictly increasing on $I_\ell-p_\ell$, which is a neighbourhood of $0$. Therefore, there exists $\delta_\ell>0$ such that 
\Eq{*}{
  \varphi\bigg(p+\sum_{j=1}^{\ell}\delta_je_j\bigg)
  -\varphi\bigg(p+\sum_{j=1}^{\ell-1}\delta_je_j \bigg) &<\frac{\varepsilon}{k}.
}
Adding up these inequalities side by side for $\ell\in\{1,\dots,k\}$, it follows that
\Eq{*}{
  \varphi\bigg(p+\sum_{j=1}^{k}\delta_je_j\bigg)
  -\varphi(p) &<\varepsilon.
}
Analogously, we can construct $\delta'_1,\dots,\delta'_k>0$ such that
\Eq{*}{
  \varphi(p)-
  \varphi\bigg(p-\sum_{j=1}^{k}\delta'_je_j\bigg)
   &<\varepsilon.
}
Therefore, since $\varphi$ is separately strictly increasing, for all $x=(x_1,\dots,x_k)\in\prod_{\ell=1}^k (p_\ell-\delta'_\ell,p_\ell+\delta_\ell)$, we have
\Eq{*}{
\varphi(p)-\varepsilon
<\varphi\bigg(p-\sum_{j=1}^{k}\delta'_je_j\bigg)
\leq\varphi(x)
\leq\varphi\bigg(p+\sum_{j=1}^{k}\delta_je_j\bigg)
<\varphi(p)+\varepsilon,
}
and hence
$$
|\varphi(x)-\varphi(p)|<\varepsilon,
$$
which shows that $\varphi$ is continuous at $p$.
\end{proof}

\Cor{Mink}{Let $k\in\N$, let $f_j:\R_+\to\R$ be strictly increasing for $j\in\{0,1,\dots,k\}$. 
Then the following conditions are equivalent.
\begin{enumerate}[(a)]
\item For all $n\in\N$ and for all $(x_{1,1},\dots,x_{1,k}),\dots,(x_{n,1},\dots,x_{n,k})\in\R_+^k$, 
\Eq{*}{
  A^{[n]}_{f_0}\big(x_{1,1}+\dots+x_{1,k},\dots,x_{n,1}+\dots+x_{n,k}\big)
  \leq A^{[n]}_{f_1}(x_{1,1},\dots,x_{n,1})+\dots+A^{[n]}_{f_k}(x_{1,k},\dots,x_{n,k}).
}
\item For all $n\in\N$, for all $(x_{1,1},\dots,x_{1,k}),\dots,(x_{n,1},\dots,x_{n,k})\in\R_+^k$, and for all $(\lambda_1,\dots,\lambda_n)\in\Lambda_n$, 
\Eq{*}{
  \widetilde{A}^{[n]}_{f_0}\big(x_{1,1}+\dots+x_{1,k}&,\dots,x_{n,1}+\dots+x_{n,k};\lambda_1,\dots,\lambda_n\big)\\
  &\leq \widetilde{A}^{[n]}_{f_1}(x_{1,1},\dots,x_{n,1};\lambda_1,\dots,\lambda_n)+\dots+\widetilde{A}^{[n]}_{f_k}(x_{1,k},\dots,x_{n,k};\lambda_1,\dots,\lambda_n).
}
\item For all $(x_1,\dots,x_k),(y_1,\dots,y_k)\in\R_+^k$, 
\Eq{*}{
  A^{[2]}_{f_0}\big(x_1+\dots+x_k,y_1+\dots+y_k\big)
  \leq A^{[2]}_{f_1}(x_1,y_1)+\dots+A^{[2]}_{f_k}(x_k,y_k).
}
\item The functions $f_0,f_1,\dots,f_k$ are continuous and the function $\Psi:f_1(\R_+)\times\dots\times f_k(\R_+)\to \R$ defined by
\Eq{*}{
  \Psi(u_1,\dots,u_k):=f_0\big(f_1^{-1}(u_1)+\dots+f_k^{-1}(u_k)\big)
}
is concave.
\end{enumerate}
}

\begin{proof} Applying \thm{Main} with 
\Eq{phi}{
\varphi(x_1,\dots,x_k):=\Phi(x_1,\dots,x_k):=x_1+\dots+x_k \qquad (x_1,\dots,x_k\in\R_+),
}
we can see that assertions $(a)$ and $(b)$
are equivalent to each other. In view of this theorem we also have that $(d)$ implies $(a)$.

By taking $\lambda_1:=\lambda_2:=1$ and $\lambda_j:=0$ for $j\geq 2$ in assertion $(b)$, we get that assertion $(c)$ holds (with an obvious substitution). 

Now assume that assertion $(c)$ is valid. Then, applying \thm{cont} with $n=2$ and $\varphi$ defined by \eq{phi}, it follows that $f_j$ is continuous for all $j\in\{0,1,\dots,k\}$. Then the generalized inverse of each of these functions is just their standard inverse and their range is convex. Using this, the inequality in $(d)$, with the substitutions $x_i=f_i^{-1}(u_i)$, $y_i=f_i^{-1}(v_i)$ yields that the function $\Psi$ defined in assertion $(d)$ is Jensen concave. Due to its continuity, it follows that it is concave. Therefore assertion $(d)$ follows.
\end{proof}

In a completely analogous manner, we can obtain the a result for Hölder-type inequalities.

\Cor{Höld}{Let $k\in\N$, let $f_j:\R_+\to\R$ be strictly increasing for $j\in\{0,1,\dots,k\}$. 
Then the following conditions are equivalent.
\begin{enumerate}[(a)]
\item For all $n\in\N$ and for all $(x_{1,1},\dots,x_{1,k}),\dots,(x_{n,1},\dots,x_{n,k})\in\R_+^k$, 
\Eq{*}{
  A^{[n]}_{f_0}\big(x_{1,1}\cdots x_{1,k},\dots,x_{n,1}\cdots x_{n,k}\big)
  \leq A^{[n]}_{f_1}(x_{1,1},\dots,x_{n,1})\cdots A^{[n]}_{f_k}(x_{1,k},\dots,x_{n,k}).
}
\item For all $n\in\N$, for all $(x_{1,1},\dots,x_{1,k}),\dots,(x_{n,1},\dots,x_{n,k})\in\R_+^k$, and for all $(\lambda_1,\dots,\lambda_n)\in\Lambda_n$, 
\Eq{*}{
  \widetilde{A}^{[n]}_{f_0}\big(x_{1,1}\cdots x_{1,k}&,\dots,x_{n,1}\cdots x_{n,k};\lambda_1,\dots,\lambda_n\big)\\
  &\leq \widetilde{A}^{[n]}_{f_1}(x_{1,1},\dots,x_{n,1};\lambda_1,\dots,\lambda_n)\cdots \widetilde{A}^{[n]}_{f_k}(x_{1,k},\dots,x_{n,k};\lambda_1,\dots,\lambda_n).
}
\item For all $(x_1,\dots,x_k),(y_1,\dots,y_k)\in\R_+^k$, 
\Eq{*}{
  A^{[2]}_{f_0}\big(x_1\cdots x_k,y_1\cdots y_k\big)
  \leq A^{[2]}_{f_1}(x_1,y_1)\cdots A^{[2]}_{f_k}(x_k,y_k).
}
\item The functions $f_0,f_1,\dots,f_k$ are continuous and the function $\Psi:f_1(\R_+)\times\dots\times f_k(\R_+)\to \R$ defined by
\Eq{*}{
  \Psi(u_1,\dots,u_k):=f_0\big(f_1^{-1}(u_1)\cdots f_k^{-1}(u_k)\big)
}
is concave.
\end{enumerate}
}

In our final corollary we characterize the Jensen convexity of generalized quasi-arithmetic means.

\Cor{Jens}{Let $f:I\to\R$ be strictly increasing. 
Then the following conditions are equivalent.
\begin{enumerate}[(a)]
\item For all $n\in\N$ and for all $(x_1,\dots,x_n),(y_1,\dots,y_n)\in I^n$, 
\Eq{*}{
  A^{[n]}_{f}\Big(\frac{x_1+y_1}{2},\dots,\frac{x_n+y_n}{2}\Big)
  \leq \frac{1}{2}\Big(A^{[n]}_{f}(x_1,\dots,x_n) +A^{[n]}_{f}(y_1,\dots,y_n)\Big).
}
\item For all $n\in\N$, for all $(x_1,\dots,x_n),(y_1,\dots,y_n)\in I^n$, and for all $(\lambda_1,\dots,\lambda_n)\in\Lambda_n$, 
\Eq{*}{
  \widetilde{A}^{[n]}_{f}&\Big(\frac{x_1+y_1}{2},\dots,\frac{x_n+y_n}{2};\lambda_1,\dots,\lambda_n\Big)\\
  &\leq \frac{1}{2}\Big(\widetilde{A}^{[n]}_{f}(x_1,\dots,x_n;\lambda_1,\dots,\lambda_n)+\widetilde{A}^{[n]}_{f}(y_1,\dots,y_n;\lambda_1,\dots,\lambda_n)\Big).
}
\item For all $(x_1,x_2),(y_1,y_2)\in I^2$, , 
\Eq{*}{
  A^{[2]}_{f}\Big(\frac{x_1+y_1}{2},\frac{x_2+y_2}{2}\Big)
  \leq \frac{1}{2}\Big(A^{[2]}_{f}(x_1,x_2)+A^{[2]}_{f}(y_1,y_2)\Big).
}
\item The function $f$ is continuous and the function $\Psi:f(I)^2\to \R$ defined by
\Eq{*}{
  \Psi(u,v):=f\bigg(\frac{f^{-1}(u)+f^{-1}(v)}{2}\bigg)
}
is concave.
\item $f$ is twice continuously differentiable with a positive first derivative and either $f''=0$ on $I$ of $f''$ is positive and $\frac{f'}{f''}$ is concave.
\end{enumerate}
}

\begin{proof}
By taking $k=2$ and functions $\varphi,\Phi:I^2\to I$ defined as
$$
\varphi(x,y):=\Phi(x,y):=\frac{x+y}{2}
$$
the equivalence of the first four assertions follows from \thm{Main} and \thm{cont} in the same manner as for the equivalence of the analogous conditions for Minkowski-type inequalities.

The equivalence of assertions $(a)$ and $(e)$ was obtained by the authors in \cite{PalPas21a}.
\end{proof}


\end{document}